\documentclass[a4paper,10pt,fleqn,enabledeprecatedfontcommands]{scrartcl}
\usepackage{amsmath}
\usepackage{amsthm}

\renewenvironment{abstract}{
  \small
  \begin{center}
    \bfseries \abstractname
  \end{center}
  \list{}{
    \setlength{\leftmargin}{7mm}
    \setlength{\rightmargin}{\leftmargin}
  }
  \item
  \relax
}{\endlist}

\usepackage{colortbl}

\usepackage{amssymb}
\usepackage{tikz}
\usetikzlibrary{patterns, spy}
\usepackage{pgfplots}

\usepackage[ruled]{algorithm2e}
\SetKwProg{Fn}{Function}{}{}
\SetKwRepeat{Do}{do}{while}

\theoremstyle{definition}

\theoremstyle{remark}

\numberwithin{equation}{section}

\newcommand{\beq}{\begin{eqnarray*}}
\newcommand{\eeq}{\end{eqnarray*}}
\newcommand{\beqn}{\begin{eqnarray}}
\newcommand{\eeqn}{\end{eqnarray}}
\newcommand{\notinclude}[1]{}

\newcommand{\R}{\mathbb{R}}
\newcommand{\N}{\mathbb{N}}

\newcommand{\domain}{D}
\newcommand{\cdomain}{\domain}

\newcommand{\Ed}{\mathbf{E}}
\newcommand{\J}{\mathcal{J}}
\newcommand{\Jd}{\mathbf{J}}

\newcommand{\Lag}{\mathbf{L}}

\renewcommand{\d}{\,\mathrm{d}}
\newcommand{\dx}{\,\mathrm{d}x}

\newcommand{\vol}{\mathcal{V}}
\newcommand{\peri}{\mathcal{L}}

\newcommand{\del}{\partial}

\renewcommand{\div}{\mathrm{div}}

\newcommand{\pf}{v}

\newcommand{\El}{\mathbf{C}}
\newcommand{\edge}{\mathbf{e}}

\newcommand{\Pf}{\mathbf{v}}

\newcommand{\eltensor}{\mathcal{C}}
\newcommand{\force}{f}

\DeclareRobustCommand\tikzsquare[1]{\tikz \draw[fill=#1] (0,0)
  rectangle (0.3,0.3);}

\definecolor{colbr1}{HTML}{00AAFF}
\definecolor{colp}{HTML}{0000DD}
\definecolor{colcp11}{HTML}{00AAAA}
\definecolor{colcp12}{HTML}{00FF88}
\definecolor{colcn1}{HTML}{00AA55}
\definecolor{colcn2}{HTML}{FF5500}
\definecolor{colcp21}{HTML}{FFAA00}
\definecolor{colcp22}{HTML}{AAAA00}
\definecolor{colbr2}{HTML}{FFDD00}
\definecolor{colcn3}{HTML}{8800FF}
\definecolor{colcp31}{HTML}{BB00FF}
\definecolor{colcp32}{HTML}{AA88AA}
\definecolor{colbr3}{HTML}{FF00FF}

\begin{document}

\newpage
\setcounter{page}{1}

\title{Branching Structures in Elastic Shape Optimization}

\author{
Nora L\"uthen \footnotemark[1]
\and Martin Rumpf \footnotemark[2]
\and Sascha T{\"o}lkes \footnotemark[2]
\and Orestis Vantzos \footnotemark[3]
}

\maketitle

\begin{abstract}
Fine scale elastic structures are widespread in nature, for instances in plants or bones, whenever stiffness and low weight are required.
These patterns frequently refine towards a Dirichlet boundary to ensure an      effective load transfer.
The paper discusses the optimization of such supporting structures
in a specific class of domain patterns in 2D, which composes of periodic and branching period transitions on subdomain facets.
These investigations can be considered as a case study to display examples of optimal branching domain patterns.\\
In explicit, a rectangular domain is decomposed into rectangular subdomains, which share facets with neighbouring subdomains or
with facets which split on one side into equally sized facets of two different subdomains.
On each subdomain one considers an elastic material phase with stiff elasticity coefficients and an approximate void phase with orders of magnitude softer material.
For given load on the outer domain boundary, which is distributed on a prescribed fine scale pattern representing 
the contact area of the shape,
the interior elastic phase is optimized with respect to the compliance cost.
The elastic stress is supposed to be continuous on the domain and a stress based finite volume discretization is used for the optimization.
If in one direction equally sized subdomains with equal adjacent subdomain topology line up, these subdomains are consider as equal copies including the enforced
boundary conditions for the stress and form a locally periodic substructure.\\
An alternating descent algorithm is employed for a discrete characteristic function describing the stiff elastic subset on the subdomains and
the solution of the elastic state equation. Numerical experiments are shown for compression and shear load on the boundary of a
quadratic domain.
  \end{abstract}

{\bfseries Key Words:} Shape optimization, elasticity, branching patterns

\renewcommand{\thefootnote}{\fnsymbol{footnote}}
\footnotetext[1]{Chair of Computational Science, Clausiusstrasse 33, ETH-Zentrum, CLT C 14, 8092 Z\"urich, Switzerland, nluethen@ethz.ch}
\footnotetext[2]{Institute for Numerical Simulation, University of Bonn, Endenicher Allee 60, 53115 Bonn, Germany, \{martin.rumpf\}\{sascha.toelkes\}@ins.uni-bonn.de}
\footnotetext[3]{Center for Graphics and Geometric Computing (CGGC), Computer Science Department, Technion, Technion City, Haifa, 32000, Israel, vantzos@cs.technion.ac.il}
\renewcommand{\thefootnote}{\arabic{footnote}}

\section{Introduction}\label{sec:introduction}
The formation of microstructures is a common phenomenon in elastic
shape optimization. We refer to \cite{Be95a} and \cite{Al02} for an
overview about these topics. Depending on the geometry of the
computational domain and the loads applied to it different type of microstructures
appear.
In fact, besides locally periodic structures one also observes branching type patterns
when optimizing with  respect to a compliance cost functional, a volume cost,
and the perimeter of the structure.
These patterns refine towards a Dirichlet boundary of the configuration \cite[Figure 13]{PeRuWi12}.
Such branching patterns can also been observed in nature, for instance in the spongiosa of bones \cite{Mu09c}.
For the basic load configurations of uniaxial load and shear load Kohn and Wirth \cite{KoWi14,KoWi15} considered scaling laws for the cost functional and for the
weight in front of the perimeter tending to zero.

Besides the question how locally periodic or branching periodic patterns might look like, a
central challenge is also to identify optimal decompositions of elastic material devices or objects
into such spatially varying patterns.
This paper should be regarded as a case study in this direction.
It is intended as a first step towards a truely multiscale modeling of optimized elastic objects
involving periodic and branching periodic patterns. Such a multiscale model would enable to apply techniques from
\emph{homogenization}, an important concept to upscale the microscopic properties
to a macroscale. This methodology has been described for instance in \cite{DoCi99, Mi02}
and in detail in the context of elastic shape optimization in \cite{Al02} and in the context of engineering applications for instance in \cite{MuBrDi15}
and \cite{NaWiAm11}.
Two-scale materials can be numerically computed by the
Heterogeneous Multiscale Method (HMM) \cite{EEn05, EEn03, EEnHu03,EMiZh05},
which explicitly simulates periodic microstructures at positions in the macroscopic domain.

Here we consider domain patterns which explicitly prescribe locally periodic or branching periodic structures on a fine scale
instead of taking into account a truely multiscale model.
To this end we study computational domains that can be
decomposed into a set of rectangular subdomains with compatibility conditions
on the facets of these domains. The actual elastic structure will be a subset of the computational domain, which we aim to identify
via elastic shape optimization. This elastic structure correspondingly splits into components on the subdomains.
Each subdomain will be a copy of one rectangular reference cell
on which the shape of a hard phase will be optimized. 
Then these elastic structures in the subdomains assemble to a locally periodic or
a locally branching type structures. The discretization discussed in this article is an
extension of the approach suggested by one of the authors, O. Vantzos and also used in \cite{Lu16}, now going
beyond purely branching periodic ensembles of cells.
\medskip

The paper is organized as follows: in Section \ref{sec:composite} we
describe the admissible subdivision of the computational domain into cells.
Then, Section \ref{sec:elastic} discusses the underlying elasticity model and the optimization problem.
Section \ref{sec:discr} presents the spatial discretization and the numerical solution of the state equation,
whereas in Section \ref{sec:pfgs}  the alternating descent scheme for the optimization of the discrete characteristic function is 
investigated.
Finally, results for two different load configurations are depicted in Section  \ref{sec:res}.

\section{Composite structures}
\label{sec:composite}

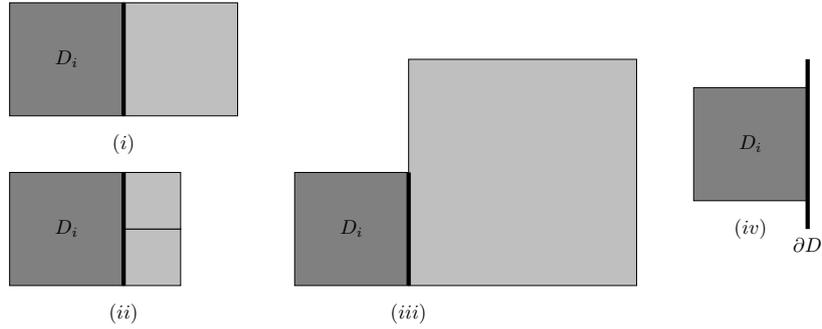
\begin{figure}[htb]
  \centering
  \scalebox{0.75}{
    \begin{tikzpicture}{x=1cm,y=1cm}
      \begin{scope}
        \draw[fill=gray] (0,0) rectangle (2,2) node [pos=0.5] {$\domain_{i}$};
        \draw[fill=lightgray] (2,0) rectangle (4,2);
        \draw[black,line width=0.75mm] (2,0) -- (2,2);
        \node at (2,-0.5) {$(i)$};
      \end{scope}
      \begin{scope}[yshift=-3cm]
        \draw[fill=gray] (0,0) rectangle (2,2)  node [pos=0.5] {$\domain_{i}$};
        \draw[fill=lightgray] (2,0) rectangle (3,1);
        \draw[fill=lightgray] (2,1) rectangle (3,2);
        \draw[black,line width=0.75mm] (2,0) -- (2,2);
        \node at (2,-0.5) {$(ii)$};
      \end{scope}
      \begin{scope}[xshift=5cm,yshift=-3cm]
        \draw[fill=gray] (0,0) rectangle (2,2) node [pos=0.5] {$\domain_{i}$};
        \draw[fill=lightgray] (2,0) rectangle (6,4);
        \draw[black,line width=0.75mm] (2,0) -- (2,2);
        \node at (2,-0.5) {$(iii)$};
      \end{scope}
      \begin{scope}[xshift=12cm,yshift=-1.5cm]
        \draw[fill=gray] (0,0) rectangle (2,2) node [pos=0.5] {$\domain_{i}$};
        \draw[black,line width=0.75mm] (2,-0.5) node[below] {$\del \cdomain$}  -- (2,2.5);
        \node at (1,-0.5) {$(iv)$};
      \end{scope}
    \end{tikzpicture}
  }
  \caption{The different local subdomain configurations at a single facet.}
  \label{fig:subdom}
\end{figure}

Let $\cdomain \in \R^2$ be a computational domain consisting of several
subdomains $\domain_{i}$, $i=1,\ldots, M$ for some $M\in \N$, such that
$\cup_{i=1}^{M} \overline{\domain_{i}} = \cdomain$ and
$\cap_{i=1}^{M} \domain_{i} = \emptyset$. Each subdomain ist supposed to be a rectangle $(a_1,a_2)\times (b_1,b_2)$.
Each facet of a subdomain
\begin{itemize}
\item[(i)]  is either also a facet of an adjacent subdomain (e.g. there is a subdomain $(a_2,a_3)\times (b_1,b_2)$
sharing the facet $\{a_2\}\times (b_1,b_2)$ with the subdomain $(a_1,a_2)\times (b_1,b_2)$\;),
\item[(ii)] or splits into two facets of two adjacent subdomains (e.g. there are subdomains $(a_2,a_3)\times (b_1,\tfrac{b_1+b_2}{2})$
and $(a_2,a_3)\times (\tfrac{b_1+b_2}{2},b_2)$ whose facets  $\{a_2\}\times (b_1,\tfrac{b_1+b_2}{2})$ and $\{a_2\}\times (\tfrac{b_1+b_2}{2},b_2)$
results from a splitting of the facet $\{a2\}\times (b_1,b_2)$\;),
\item[(iii)] or is on the facets resulting from such a splitting of a facet of an adjacent subdomain,
\item[(iv)] or is a boundary facet.
\end{itemize}
We always  assume that the splitting of facets is in two halves of
equal length. The subdomain configurations at a single facet are shows in Figure~\ref{fig:subdom}.
The four facets of a single rectangular subdomain can be of different type.

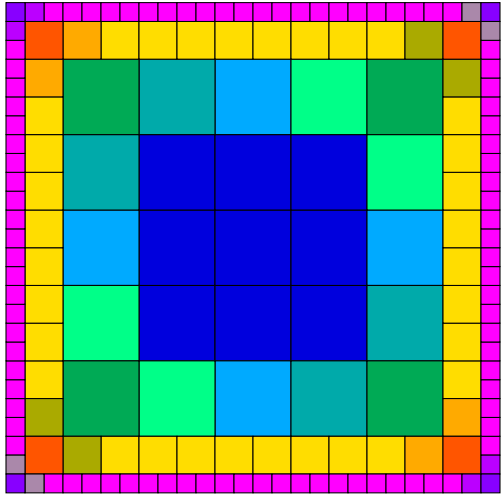
\begin{figure}[htb]
  \centering
  \begin{tikzpicture}
    \foreach \i in {0,1,2} {
      \foreach \j in {0,1,2} {
        \begin{scope}[shift={(\i,\j)}]
          \draw[fill=colp] (0,0) rectangle (1,1);
        \end{scope}
      }
    }
    \foreach \i in {-1,3} {
      \begin{scope}[shift={(\i,1)}]
        \draw[fill=colbr1] (0,0) rectangle (1,1);
      \end{scope}
      \begin{scope}[shift={(1,\i)}]
        \draw[fill=colbr1] (0,0) rectangle (1,1);
      \end{scope}
    }
    \foreach \i/\j in {{0}/{3},{-1}/{2},{2}/{-1},{3}/{0}} {
      \begin{scope}[shift={(\i,\j)}]
        \draw[fill=colcp11] (0,0) rectangle (1,1);
      \end{scope}
    }
    \foreach \i/\j in {{2}/{3},{-1}/{0},{0}/{-1},{3}/{2}} {
      \begin{scope}[shift={(\i,\j)}]
        \draw[fill=colcp12] (0,0) rectangle (1,1);
      \end{scope}
    }
    \foreach \i in {-1,3} {
      \foreach \j in {-1,3} {
        \begin{scope}[shift={(\i,\j)}]
          \draw[fill=colcn1] (0,0) rectangle (1,1);
        \end{scope}
      }
    }
        \foreach \i in {-0.5,0,...,3} {
      \foreach \j in {4, -1.5} {
        \begin{scope}[shift={(\i,\j)},scale=0.5]
          \draw[fill=colbr2] (0,0) rectangle (1,1);
        \end{scope}
        \begin{scope}[shift={(\j,\i)},scale=0.5]
          \draw[fill=colbr2] (0,0) rectangle (1,1);
        \end{scope}
      }
    }
    \foreach \i/\j in {{-1}/{4},{-1.5}/{3.5},{3.5}/{-1.5},{4}/{-1}} {
      \begin{scope}[shift={(\i,\j)},scale=0.5]
        \draw[fill=colcp21] (0,0) rectangle (1,1);
      \end{scope}
    }
    \foreach \i/\j in {{3.5}/{4},{-1.5}/{-1},{-1}/{-1.5},{4}/{3.5}} {
      \begin{scope}[shift={(\i,\j)},scale=0.5]
        \draw[fill=colcp22] (0,0) rectangle (1,1);
      \end{scope}
    }
    \foreach \i/\j in {{-1.5}/{4},{4}/{4},{-1.5}/{-1.5},{4}/{-1.5}} {
      \begin{scope}[shift={(\i,\j)},scale=0.5]
        \draw[fill=colcn2] (0,0) rectangle (1,1);
      \end{scope}
    }
        \foreach \i in {-1.25,-1,...,4} {
      \foreach \j in {4.5, -1.75} {
        \begin{scope}[shift={(\i,\j)},scale=0.25]
          \draw[fill=colbr3] (0,0) rectangle (1,1);
        \end{scope}
        \begin{scope}[shift={(\j,\i)},scale=0.25]
          \draw[fill=colbr3] (0,0) rectangle (1,1);
        \end{scope}
      }
    }
    \foreach \i/\j in {{-1.5}/{4.5},{-1.75}/{4.25},{4.25}/{-1.75},{4.5}/{-1.5}} {
      \begin{scope}[shift={(\i,\j)},scale=0.25]
        \draw[fill=colcp31] (0,0) rectangle (1,1);
      \end{scope}
    }
    \foreach \i/\j in {{4.25}/{4.5},{-1.75}/{-1.5},{-1.5}/{-1.75},{4.5}/{4.25}} {
      \begin{scope}[shift={(\i,\j)},scale=0.25]
        \draw[fill=colcp32] (0,0) rectangle (1,1);
      \end{scope}
    }
    \foreach \i/\j in {{-1.75}/{4.5},{4.5}/{4.5},{-1.75}/{-1.75},{4.5}/{-1.75}} {
      \begin{scope}[shift={(\i,\j)},scale=0.25]
        \draw[fill=colcn3] (0,0) rectangle (1,1);
      \end{scope}
    }
  \end{tikzpicture}
  \caption[Sketch of a domain decomposite]{
An examplary decomposition of the computational domain is displayed
with 13 different reference domains that are mapped to the subdomains of the decomposition,
where identical colors classify the different reference domains.
The different types are locally periodic cells
    (\tikzsquare{colp}), branching  cells (\tikzsquare{colbr1}
    \tikzsquare{colbr2} \tikzsquare{colbr3}), double branching 
    corner cells (\tikzsquare{colcn1} \tikzsquare{colcn2}
    \tikzsquare{colcn3}) and coupling cells (\tikzsquare{colcp11}
    \tikzsquare{colcp12} \tikzsquare{colcp21} \tikzsquare{colcp22}
    \tikzsquare{colcp31} \tikzsquare{colcp32}).}
  \label{fig:csketchcomp}
\end{figure}

Let us assume that each subdomain $\domain_{i}$ contains
an inscribed subpart of the actual elastic object, the shape of which will be optimized.
Let us denote by $\chi_i$ the associated characteristic function.
for which we consider the continuous extension in $BV$ onto the closure of the subdomain.
Thus, $\sum_{i=1}^M \chi_i$ as function on $\overline{\cdomain}$
is the characteristic function of the elastic objects we are investigating.
Each (geometric) subdomain $\domain_{i}$ will have a reference
domain assigned and several geometric subdomains can have the same
reference domain. The referenced domain is mapped onto a geometric subdomain via a translation
and a rotation by a multiple of $\tfrac{\pi}{2}$. The characteristic function $\chi_i$ and later also the force distribution on $[\chi_i=1]$ is
handled and updated on the associated reference domain.
Thus, the computational complexities scales with the number of reference domains.
If adjacent subdomains share the same reference domain and the same rotation, then they are building blocks of local, period elastic
structures.
Fig.~\ref{fig:csketchcomp} illustrates such a domain decomposition into rectangular boxes.

As a consequence, the domain sketched in
Figure~\ref{fig:csketchcomp} can later be numerically optimized by a shape optimization on
13 (coupled) reference domains only.

\section{Elastic state equation and compliance optimization}
\label{sec:elastic}

To solve the elastic state equation in the context of the shape optimization
we use the quadratic energy of a stress based formulation of linearized elasticity.
On the same basis we evaluate the compliance-type target functional.
On each reference domain the shape will be
modelled by a phase field $\pf:\domain_{i} \to \R$ where $\pf$ is assumed to approximate the characteristic function of the elastic object $\chi$.
Thus, $\{x\in \domain_{i} \,|\, \pf(x) \approx 1\}$ corresponds to the actual elastic domain
and $\{x\in \domain_{i} \,|\, \pf(x) \approx 0\}$ to the void (or in our calculation very soft) phase.

Let us consider elastic stresses $\tau: \domain_{i} \to \R^{2}$ as extensions of the elastic stresses on the elastic object $[\chi=1]$.
Let us recall that stresses act on normals of infinitesimal area elements and represent the force density acting on this area element.
We apply solely boundary forces and no volume forces. Hence, elastic stresses are divergence free, which constitutes together with the boundary condition
the state equation of our optimization problem.
Given the elasticity tensor $C(\pf)$, which depends on the phase (hard or soft) the stored elastic energy elastic on $\domain_{i}$ is defined as
\begin{equation}
  \label{eq:strelstr}
  \mathcal{E}[\pf,\tau] = \int_{\domain_{i}} \eltensor^{-1}(\pf) \tau : \tau \dx\,.
\end{equation}
Here $\eltensor(\pf)$ is a fourth order tensor satisfying $\eltensor_{ijkl}(\pf) = \eltensor_{jikl}(\pf) = \eltensor_{ijlk}(\pf) = \eltensor_{klij}(\pf)$.
In the context of this paper, we define
\[
\eltensor(\pf) = v \eltensor_{NL}
\]
where $\eltensor_{NL}$ is the elasticity tensor of the linearized
Navier-Lam{\'e} elasticity model and $\delta > 0$. We define the inverse
$\eltensor^{-1}$ (needed in \eqref{eq:strelstr}) using
\emph{Young's modulus} $E$ and the Poisson ratio $\nu$ (in Voigt
notation)
\begin{equation*}
  \eltensor^{-1}_{NL} = \frac{1}{E}
  \begin{pmatrix}
    1 & -\nu & 0\\
    -\nu & 1 & 0\\
    0 & 0 & 2 + 2\nu
  \end{pmatrix}\,.
\end{equation*}
For the phase field function $\pf$ we consider a Modica Mortola type functional (cf. \cite{MoMo77})
\begin{equation*}
  \peri^\varepsilon[\pf] = \int_{\domain_{i}} \frac{1}{\varepsilon} W(\pf) +
  \frac{\varepsilon}{2} \left| \bigtriangledown \pf \right|^2 \dx
\end{equation*}
which approximates the length of the interface between hard and soft material.
Here,
\begin{equation*}
  W(v) = \begin{cases}
    \frac{32}{\pi^2} (1-v)(v-\delta) & v \in [\delta, 1]\\
    \infty & \text{else}
  \end{cases}
\end{equation*}
denotes a double-well potential, i.\,e. a positive function
that is attending its only two minima at the pure phases $\pf = \delta$ and $\pf = 1$.  We assume that $v|_{\partial \cdomain}$ is prescribed and describes the imposed
fine scale structure on the domain boundary. The parameter $\epsilon$ is proportional to the width of the diffused interface.
Futhermore, given the phase field $\pf$ we can easily compute an approximation of the area of the elastic object
$\vol^\varepsilon[\pf] = \int_{\domain_{i}} \pf(x) \dx$.

Combining the energies above, we obtain the objective functional
\begin{equation}\label{eq:objfctnl}
  \J[\pf] = \min_{\tau \in \Sigma_\text{ad}} \mathcal{E}[\pf,\tau] + \beta \vol^\varepsilon[\pf] + \eta \peri^\varepsilon[\pf]\,,
\end{equation}
of our constraint optimization problem,
where the stress $\tau$ is minized over the set of admissable stresses
\begin{equation}
  \label{eq:defsigmaad}
  \Sigma_{\text{ad}} = \{ \sigma: \domain_{i} \to \R^{2,2}\,|\, \sigma = \sigma^T, \div\, \sigma = 0, \text{b.c.} \}\,.
\end{equation}
Here, the boundary condition (b.c.) differs for facets on $\partial \cdomain$ and for interior facets, where the stresses are continuous  across the facet.
For facets on $\partial \cdomain$  we prescribe forces $f$ with $\tau(x) \cdot n(x) = f(x)$
if $\pf(x) =1$. Here, $n(x)$ denotes the outer normal at points $x\in \partial \cdomain$.

As in \cite{KoWi14,KoWi15}, we consider the case of a vanishing Poisson
ratio $\nu = \frac{\lambda}{2(\lambda + \mu)} = 0$. As expressed by Kohn and Wirth in \cite{KoWi14}, this restriction is not expected to have a strong influence because for
truss-like structures the lateral contraction is less relevant.
The tensor $\eltensor^{-1}$ then reduces to
\begin{equation*}
  \eltensor^{-1}_{NL} = \frac{1}{E}
  \begin{pmatrix}
    1 & 0 & 0\\
    0 & 1 & 0\\
    0 & 0 & 2
  \end{pmatrix}\,.
\end{equation*}
Thus, we obtain for the stored elastic energy on the subdomain $\domain_i$
\begin{equation*}
 \mathcal{E}[\pf,\tau]   =
  \int_{\domain_i} \frac1{\pf \,E} (\tau_{11}^2 + \tau_{22}^2) +
  \frac{2}{\pf \,E} \tau_{12}^2 \dx = \int_{\domain_i} \frac{1}{\pf \,E} |\tau|^2 \dx\,.
\end{equation*}
In what follows we choose $E=1$.

\section{Discretization}
\label{sec:discr}

In this section, we present a finite volume discretization of the state equation and discuss the
optimization algorithm. For details we refer to \cite{Lu16}.
To this end, the conditions for $\sigma$ prescribed in the set of admissible stress fields
$\Sigma_{\text{ad}}$ in \eqref{eq:defsigmaad} have to be transcribed into a linear system of equations
resulting from the finite volume discretization.
For a subdomain  $\domain_{i}$ we consider a finite decomposition of $\domain_{i}$
into $N \times N$ rectangular cells $\El$ of equal size.  
The discrete forces are defined as constant vectors in $\R^2$ on edges and
the discrete phase field is assumed to be constant on cells.
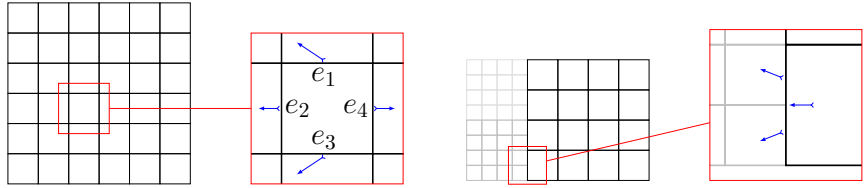
\begin{figure}[htb]
  \centering
  \scalebox{0.4}{
    \begin{tikzpicture}[spy using outlines={rectangle, red,
        magnification=3, size=5cm, connect spies}]{x=1cm,y=1cm}
      \foreach \x in {0,...,5} {
        \foreach \y in {0,...,5} {
          \draw (\x, \y) rectangle (\x + 1, \y + 1);
        }
      }

      \spy on (2.5,2.5) in node [left] at (13,2.5);
      \begin{scope}[xshift=8cm,yshift=2.5cm,scale=3]
        \draw[>-{latex},blue,line width=0.3mm] (0.8,0.53) -- (0.5,0.73);
        \draw[>-{latex},blue,line width=0.3mm] (0.31,0) -- (0.08,0);
        \draw[>-{latex},blue,line width=0.3mm] (0.8,-0.53) -- (0.5,-0.73);
        \draw[>-{latex},blue,line width=0.3mm] (1.35,0) -- (1.58,0);

        \node[scale=2] at (0.8,0.35) {\Large $e_1$};
        \node[scale=2] at (0.51,0) {\Large  $e_2$};
        \node[scale=2] at (0.8,-0.35) {\Large  $e_3$};
        \node[scale=2] at (1.15,0) {\Large  $e_4$};
      \end{scope}
    \end{tikzpicture}
  }
  \hspace{0.5cm}
  \scalebox{0.4}{
    \begin{tikzpicture}[spy using outlines={rectangle, red,
        magnification=4, size=5cm, connect spies}]{x=1cm,y=1cm}
      \foreach \x in {0,0.5,...,1.5} {
        \foreach \y in {0,0.5,...,1.5} {
          \draw[lightgray] (\x, \y) rectangle (\x + 0.5, \y + 0.5);
          \draw[lightgray!50] (\x, \y + 2) rectangle (\x + 0.5, \y + 2.5);
        }
      }
      \foreach \x in {2,...,5} {
        \foreach \y in {0,...,3} {
          \draw (\x, \y) rectangle (\x + 1, \y + 1);
        }
      }

      \spy on (2,0.5) in node [left] at (13,2.5);

      \begin{scope}[xshift=7cm,yshift=2.5cm,scale=3]
        \draw[>-{latex},blue,line width=0.3mm] (1.13,0.3) -- (0.88,0.4);
        \draw[>-{latex},blue,line width=0.3mm] (1.13,-0.3) -- (0.88, -0.4);
        \draw[>-{latex},blue,line width=0.3mm] (1.46926, 0) -- (1.2, 0);
      \end{scope}
    \end{tikzpicture}
  }
  \caption{Discretization of forces as vector valued degrees of freedom
  on the interior cells of the finite volume grid (left) and on a facet
  with prescribed branching-type boundary conditions (right). The numbering of the edges is displayed on the left.}
  \label{fig:discforcesinner}
\end{figure}
As degrees of freedom we consider average forces across edges $f_j$ which approximate
$\frac1{h_j} \int_{\edge_j} \sigma \cdot n_j \d l$ of a volume cell $\El$ with edge length $h_j$.
(cf. Fig.~\ref{fig:discforcesinner}  for a sketch).
The symmetry of the stress tensor is transformed into a conservation of torque:
\begin{align*}
  0 &= \int_{\partial \El} x \times (n \cdot \sigma) \d l=
  \int_{\partial \El} n \cdot (\sigma x^{\bot}) \d l = \int_\El
      \div(\sigma x^\bot) \dx\\
  = &\int_\El \div \sigma \,x^\bot + \sigma^T : \bigtriangledown x^\bot
  \dx = \int_\El \sigma_{21} - \sigma_{12}\dx\,,
\end{align*}
with $x= (x_1, x_2)$ and $x^\bot = (-x_2, x_1)$ for $\div\, \sigma=0$.
The conservation law $\div\, \tau =0$ in the continuous set up translates to a balance relation on cells, i.e. one obtains
\begin{align*}
  0 &= \int_\El \div \tau \d x = \int_{\partial \El} n \cdot
      \tau \d l = \sum_{j = 1}^4 \int_{\edge_j} n \cdot \tau \d l
  \approx \sum_{i=1}^4 h_i f_i \,
\end{align*}
where the discrete forces $f_i$ for $i=1,\ldots, 4$ are associated with the four edges $e_1,\ldots, e_4$ (numbered counter clock wise, starting with the upper edge, cf. Figure~\ref{fig:discforcesinner}) of the cell $\El$ and $h_i$ are the corresponding edge lengths.
Thus, the discrete balance of forces reads as
\begin{align} \label{eq:fb}
0= \sum_{i=1}^4 h_i f_i
\end{align}
for all cells $\El$. Given the above numbering of the edges, the balance of torques turns into the equation
\begin{align} \label{eq:torques}
0= (f_{1,1}+f_{4,1})- (f_{2,2}+f_{3,2})
\end{align}
for all cells $\El$, where $\force_{i,j}$ is the $j$th component of the force vector $\force_i$.
Finally, the discrete boundary conditions are encoded as follows. We consider a piecewise constant force $f$ for facets of subdomains on the boundary $\partial \cdomain$.
These forces are then equally distributed on the edges of cells $\El$ touching $\partial \cdomain$ on which $v=1$. In explicit,
given a subdomain $\domain_i$ with a facet $F$ on $\partial \cdomain$  and an element $E$ with $v=1$ and an edge $e_j$ on $F$, we define the force density
\begin{align} \label{eq:fbc}
f_j = \frac{\int_{F} v \d l }{\int_{F}  \d l} f
\end{align}
All the discrete counterparts of the conditions in \eqref{eq:defsigmaad} are assembled in a linear system $\mathbf{A} \mathbf{f} = \mathbf{b}$ for the vector of forces $\mathbf{f}$ on all edges. The matrix $\mathbf{A}$ and the right hand side $\mathbf{b}$ have the following block structure
\begin{align*}
  \mathbf{A}  =
  \begin{bmatrix}
   \mathbf{A}_{f} \\    \mathbf{A}_{t} \\   \mathbf{A}_{bc} \\
  \end{bmatrix}\,,\quad
    \mathbf{b}  =
  \begin{bmatrix}
   0 \\   0 \\   \mathbf{b}_{bc} \\
  \end{bmatrix}\,.
\end{align*}
Here, the index $f$ refers to the force balance, the index $t$ to the conservation of torque, and $bc$ to the boundary condition.
Now, solving the state equation coincides with minimizing the stored elastic energy
\begin{align*}
\Ed = \frac12 \sum_\El h^2 \sum_{i=1}^4 \sum_{j=1}^2 \frac{(f_{i,j}(\El))^2}{E}
 \end{align*}
 with $\force_{i}(\El)$ denoting the force vectors on the edges of the cell $\El$,
subject ot the constraint $\mathbf{A} \mathbf{f} = \mathbf{b}$.
This can be rephrased in a Lagrangian formulation for the Lagrangian
\begin{equation*}
  \Lag(\mathbf{f},\mathbf{\lambda}) = \frac12 \mathbf{f}^T \mathbf{M} \mathbf{f} + \mathbf{f} \mathbf{\lambda}^T(\mathbf{b} - \mathbf{A}\mathbf{f})\,.
\end{equation*}
Due to the invertibility of the matrix $\mathbf{M}$ we obtain the equations
\begin{align}
  \mathbf{f} &= \mathbf{M}^{-1} \mathbf{A}^T \mathbf{\lambda} \quad \text{ and} \nonumber\\
  \label{eq:forceseqsys}
  \mathbf{A}\mathbf{M}^{-1}\mathbf{A}^T\mathbf{\lambda} &= \mathbf{b}\,.
\end{align}
as the necessary conditions for a saddle point. They have to be solved first for the dual solution $\mathbf{\lambda}$ and then for the force vector $\mathbf{f}$.
If  $\mathrm{ker} \mathbf{A}^T=\{0\}$, then $\mathbf{Z} =  \mathbf{A}\mathbf{M}^{-1}\mathbf{A}^T$ is positive definite and thus invertible.
In general $\mathbf{A}^T$ is underdetermined with nontrivial kernel. Thus, we have to eliminate rows of $\mathbf{A}$ to reduce the kernel.
For domains consisting of a single subdomain $\domain_i$ with either
periodic, branching periodic or non-periodic boundary conditions, the
linear dependencies can be reduced to a small number of cases
(cf. \cite{Lu16}). For composite domains consisting of several subdomains
with different coupling and boundary conditions, the
number, type and complexity of linear dependencies increases.
Let us suppose that the number of boundary conditions (e.g. the number of edges on the boundary with prescribed forces is large enough to ensure that
$\mathbf{A}^T$ has more rows than columns.
Then, we apply a  QR-decomposition of $\mathbf{A}^T$
to find a basis of its kernel and use this to reduce the linear system $\mathbf{A} \mathbf{f} = \mathbf{b}$
by the elimination of redundant equations. This in turn successively reduces
$\mathbf{A}$ and $\mathbf{b}$ until $A$ has full (column-)rank and thus $\mathrm{ker} \mathbf{A}^T=\{0\}$.
For the results presented here the \texttt{CHOLMOD}
\cite{ChDaHa08} as a tool in the \texttt{Suitesparse}
package was used, which comes along with an efficient parallel implementation.

\section{Shape optimization}
\label{sec:pfgs}
Optimization of the elastic shape coincides in the discrete set up with an optimization of the discrete phase field and thus
a minimization of the objective functional \eqref{eq:objfctnl} added up over all subdomains. 
To this end, we apply an alternating solution strategy, i.e.
we alternatingly solve for the forces $\mathbf{f}$ for fixed discrete phase field $\Pf$ and improve the phase field $v$ for given forces $\mathbf{f}$.
A threshold for the difference between two consecutive phase fields in $L^2$ is taken into account as a stopping criterium for this descent scheme.
To improve the discrete phase field given a force vector $\mathbf{f}$ we apply a
Gauss-Seidel type iteration that optimizes the values of $\Pf$ on single cells of the finite volume mesh.
Let us consider the discretized version of $\J$ given in \eqref{eq:objfctnl} 
\begin{align*}
  \Jd[\Pf] := &\min_{\tau \in \Sigma_\text{ad}} \sum_\El  \Jd_\El[\Pf]\; \text{with}\\
&  \Jd_\El[\Pf] = 
    \frac{|\sigma(\El)|^2}{\Pf(\El)} + \beta \Pf(\El) + \frac{\eta}
    \epsilon \frac{32}{\pi^2} (\Pf(\El) - \delta)(1 - \Pf(\El))\\
  &\qquad \quad \;+ \frac{\eta \epsilon}{4} \sum_{i = 1}^4 \frac{(\Pf(\El) -
    \Pf(\El^{(i)}))^2}{h^2}
\end{align*}
Here, $\El^{(i)}$ denotes the cell adjacent to $\El$ across the edge $\edge_i$.
To find a minimum of $\Jd[\Pf]$ for $\Pf(\El)$ for all $\El$ and for fixed $\Pf(\El^{(i)})$, we
we apply Newton's method to compute the minimum of the rescaled function 
$\Pf(\El) \mapsto \Jd_\El(e^{-\Pf})$. The rescaling turns out to be an appropriate reformulation 
to overcome difficulties due to the singularity at $v(\El) = 0$.

To ensure a good performance of Newton's method, it is important to
choose a suitable initialization of $v$. To this end we 
consider the local terms of the cost functional dropping the term involving the discrete gradient and define:
$\tilde \Jd_\El(\Pf) = \frac{|\sigma|^2}{\Pf} + \beta \Pf + \frac\eta\epsilon \frac{32}{\pi^2} (\Pf - \delta)(1 - \Pf)\,$.
It is easy to check that this function has two minima and we initialize $\Pf$ 
with the minimal value of these two minima, again using in the algorithm the rescaled function $\Pf(\El) \mapsto \tilde \Jd_\El(e^{-\Pf})$ (or an approximation of it)
as long as this value is smaller or equal $1$. 
Otherwise, we initialize with the value $1$. 
This led to an initialization $\Pf$ depending on the current stress value on $\El$.
Starting from the second iteration of the alternating descent
algorithm, the phase field could also be initialized using the result of
the last phase field optimization. However, this approach is prone to
becoming stuck in local minima. 

\section{Results}
\label{sec:res}
\begin{figure}
  \centering
  \begin{tikzpicture}
    \draw[lightgray,fill=lightgray] (1,1) rectangle (4,4);
    \foreach \i in {0,...,5} {
      \draw[->] (1 + \i/5*3, 4.6) -- (1 + \i/5*3, 4.1);
      \draw[->] (0.4,1 + \i/5*3) -- (0.9, 1 + \i/5*3);
      \draw[->] (1 + \i/5*3, 0.4) -- (1 + \i/5*3, 0.9);
      \draw[->] (4.6,1 + \i/5*3) -- (4.1, 1 + \i/5*3);
    }

    \draw[lightgray,fill=lightgray] (7,1) rectangle (10,4);
    \foreach \i in {0,...,4} {
      \draw[->] (7 + \i/5*3, 4.2) -- (7.5 + \i/5*3, 4.2);
      \draw[->] (6.8,1.5 + \i/5*3) -- (6.8, 1 + \i/5*3);
            \draw[->] (7.5 + \i/5*3, 0.8) -- (7 + \i/5*3, 0.8);
      \draw[->] (10.2, 1 + \i/5*3) -- (10.2, 1.5 + \i/5*3);
    }
  \end{tikzpicture}
  \caption{Sketch of the applied loads: Compression (left) and shear (right)}
  \label{fig:comploads}
\end{figure}
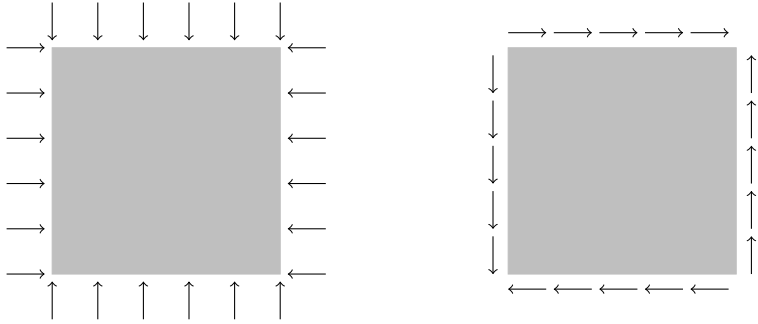
We applied our method for two different loads on a quadratic domain
$\domain \subset [0, 1]^2$, subdivided into subdomain as displayed in Figure~\ref{fig:csketchcomp}.
In particular, the central region of the domain is filled with  subdomains with inscribed locally periodic elastic structure and periodic boundary conditions for the forces,
where as in the vicinity of the boundary branching periodic subdomains are taken into account.
The two load scenarios are compression and shear as depicted in  Figure~\ref{fig:comploads}.

For all domain types, the subdomains were discretized using $N \times
N$ cells with $N = 200$. On each cell of the
outer layer, forces were applied on the edge intervals $[\frac26 N, \frac36
N]$, $[\frac46 N, \frac56 N]$ on all horizontal and vertical boundaries.
These forces are depicted in
Figure~\ref{fig:comploads} for both scenario.

\begin{figure}[h]
  \centering
  \scalebox{0.75}{
    \begin{tikzpicture}[spy using outlines={rectangle, black,
        magnification=3, size=3cm, connect spies}]
      \node{\pgfimage[height=10cm]{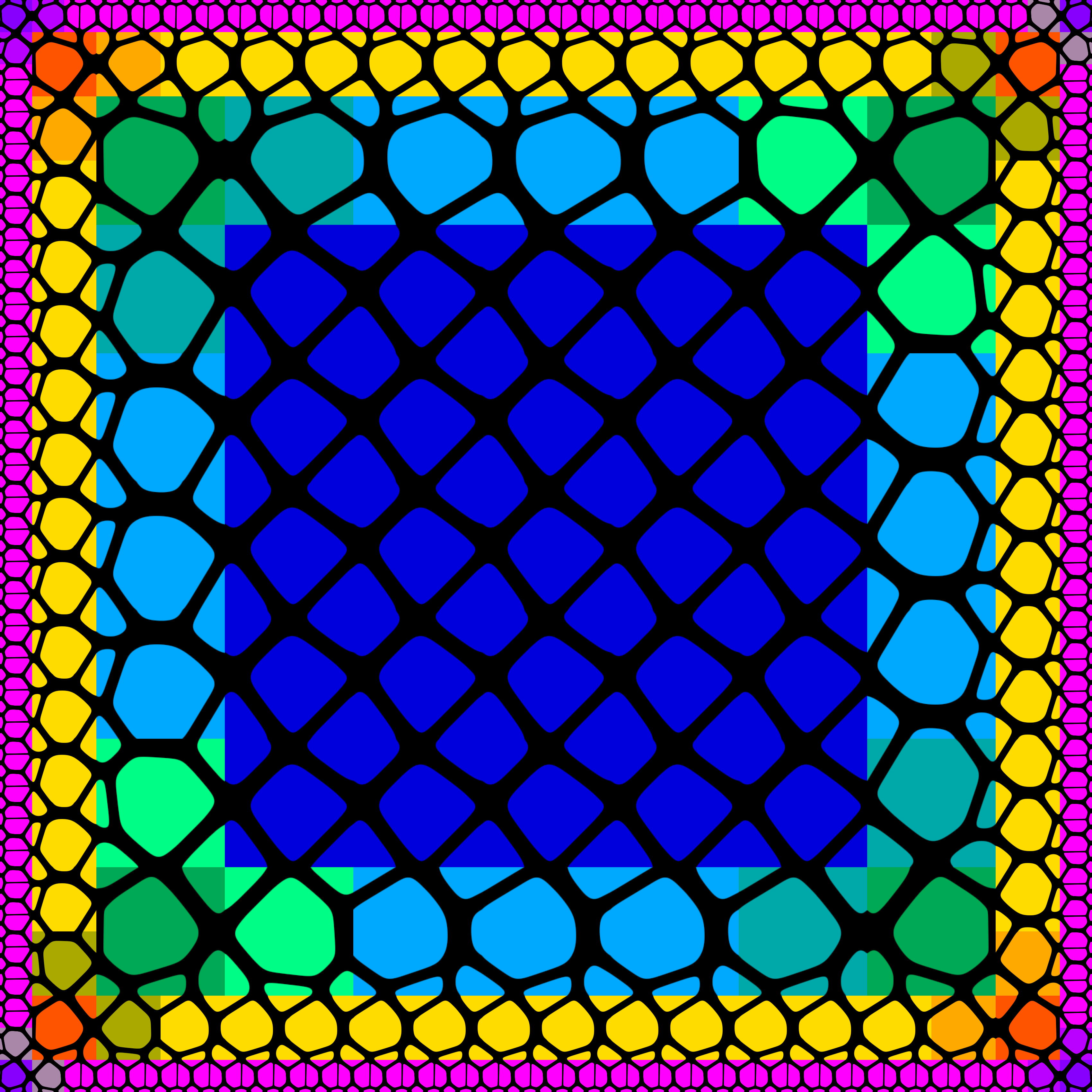}};
      \spy on (4.5,4.5) in node [left] at (10,4);
      \spy on (4.5,0) in node [left] at (10,0);
      \spy on (2.9, -3.6) in node [left] at (10, -4);
    \end{tikzpicture}
  }
  \vspace{0.5cm}

  \scalebox{0.75}{
    \begin{tikzpicture}[spy using outlines={rectangle, black,
        magnification=3, size=3cm, connect spies}]
      \node{\pgfimage[height=10cm]{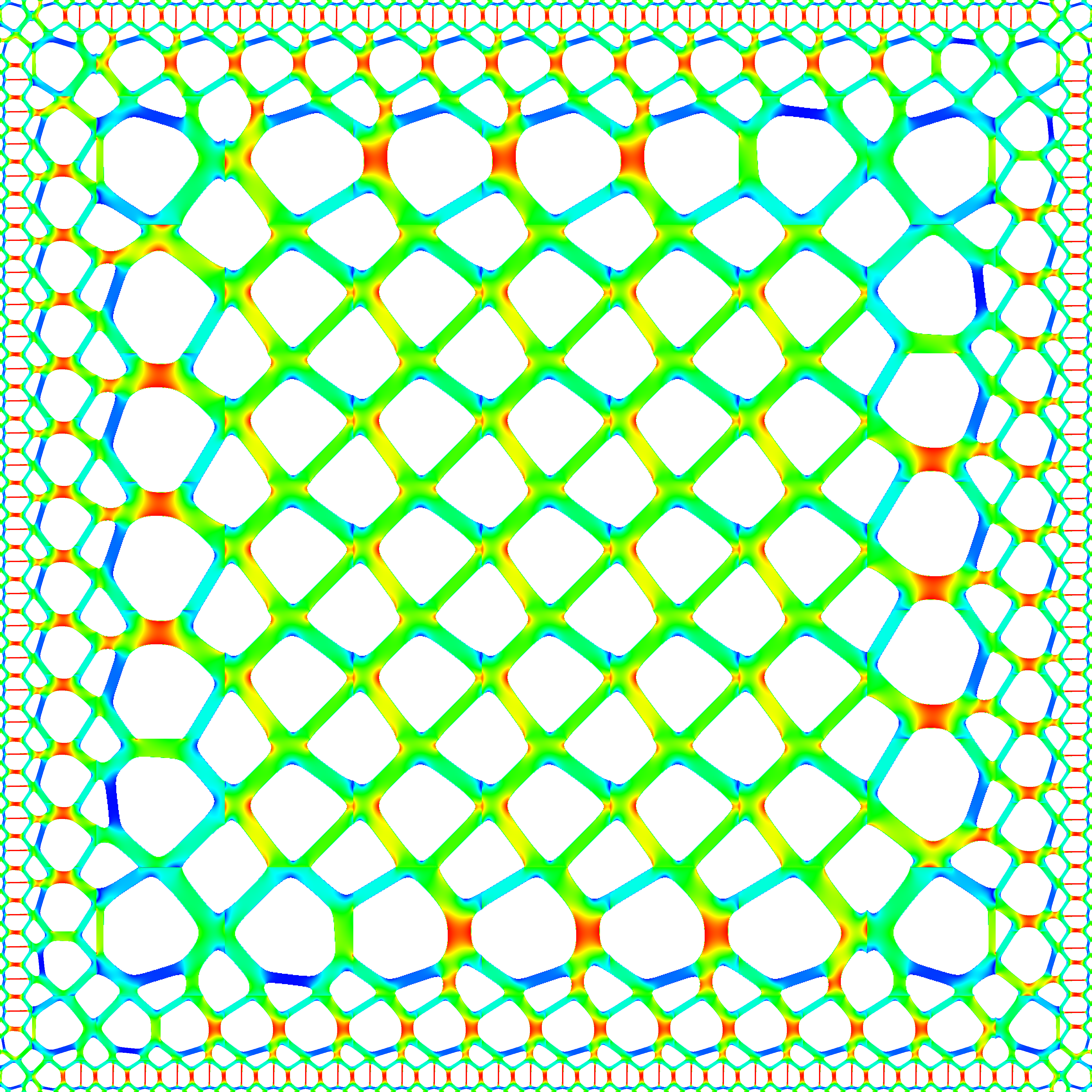}};
      \spy on (4.5,4.5) in node [left] at (10,4);
      \spy on (4.5,0) in node [left] at (10,0);
      \spy on (2.9, -3.6) in node [left] at (10, -4);
    \end{tikzpicture}
  }
  \caption{The optimal shape for the compression load case and a
    subdomain structure with $13$ cell types is depicted (top).
    We show in addition a color coding of the von Mises stresses using the
    colorbar \protect\includegraphics[height=1ex,width=5em]{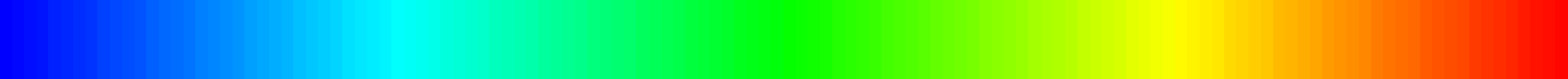}
    on the set $[\Pf > 0.5]$ (bottom).}
  \label{fig:compresc13}
\end{figure}

{\bfseries Compression load:}
Pillar-like structures support the load on the boundary,
branching structures transfer load to a mesh-like structure
in the center region of the domain as depicted in Fig~\ref{fig:compresc13}. 
The elastic structures on corner and coupling cells connect the
branching periodic pillars. Let us point to a small artifact in the optimal shape shown in one of the magnifications, 
where the phase field could obviously not been fully optimized locally.

\begin{figure}[h]
  \centering
  \scalebox{0.75}{
    \begin{tikzpicture}[spy using outlines={rectangle, black,
        magnification=3, size=3cm, connect spies}]
      \node{\pgfimage[height=10cm]{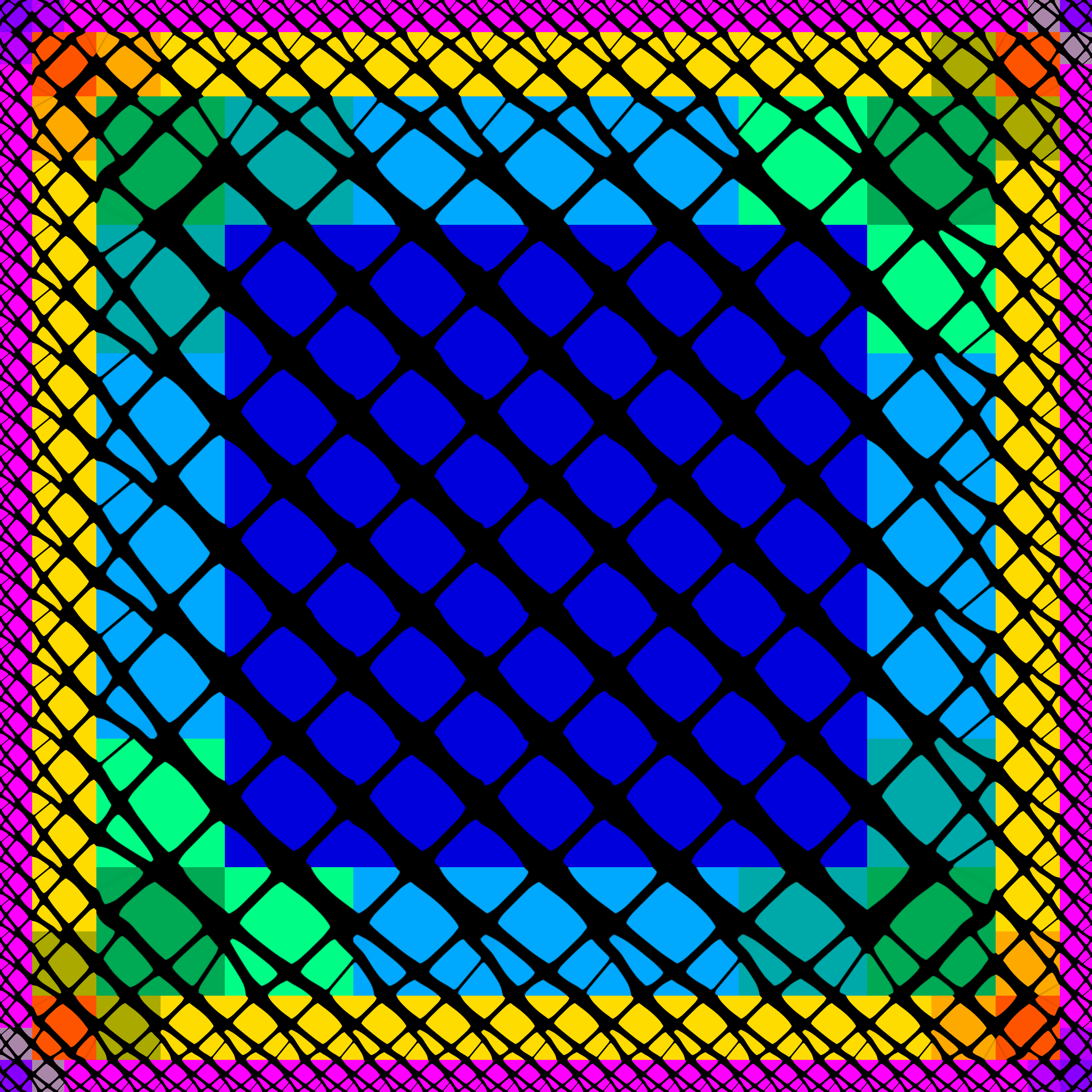}};
      \spy on (4.5,4.5) in node [left] at (10,4);
      \spy on (4.5,0) in node [left] at (10,0);
    \end{tikzpicture}
  }\\[0.5cm]
  \scalebox{0.75}{
    \begin{tikzpicture}[spy using outlines={rectangle, black,
        magnification=3, size=3cm, connect spies}]
      \node{\pgfimage[height=10cm]{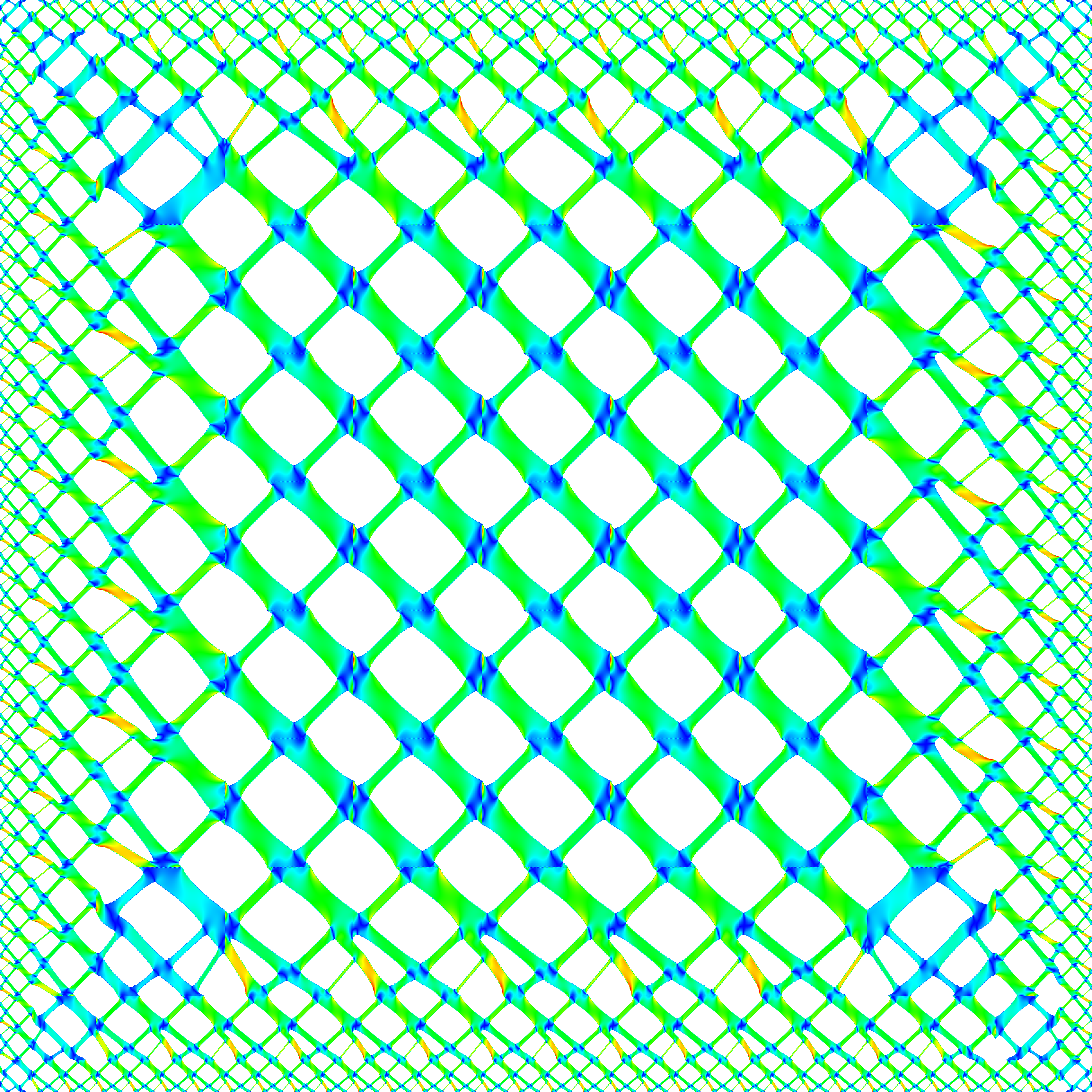}};
      \spy on (4.5,4.5) in node [left] at (10,4);
      \spy on (4.5,0) in node [left] at (10,0);
    \end{tikzpicture}
  }
  \caption{The optimal shape for the shear load case and the same
    subdomain structure as in Fig.~\ref{fig:compresc13} is displayed. 
 with two regions being magnified on the right. Again together with
 the optimal shape (top) the associated von Mises stresses are
 rendered (bottom).}
  \label{fig:shear13}
\end{figure}

{\bfseries Shear load:} Fig~\ref{fig:shear13} shows the optimal elastic shape in case of the shear load.
On the resulting compound, optimized elastic structure forces are transferred via branching type structures to the center.

\renewcommand{\refname}{References}
\bibliographystyle{acm}
\bibliography{bibtex/all,bibtex/own}

\def\polhk#1{\setbox0=\hbox{#1}{\ooalign{\hidewidth
  \lower1.5ex\hbox{`}\hidewidth\crcr\unhbox0}}}
\begin{thebibliography}{10}

\bibitem{Al02}
{\sc Allaire, G.}
\newblock {\em Shape optimization by the homogenization method}, vol.~146 of
  {\em Applied Mathematical Sciences}.
\newblock Springer-Verlag, New York, 2002.

\bibitem{Be95a}
{\sc Bends{\o}e, M.~P.}
\newblock {\em Optimization of structural topology, shape, and material}.
\newblock Springer-Verlag, Berlin, 1995.

\bibitem{ChDaHa08}
{\sc Chen, Y., Davis, T.~A., Hager, W.~W., and Rajamanickam, S.}
\newblock Algorithm 887: Cholmod, supernodal sparse cholesky factorization and
  update/downdate.
\newblock {\em ACM Trans. Math. Softw. 35}, 3 (Oct. 2008), 22:1--22:14.

\bibitem{DoCi99}
{\sc Cioranescu, D., and Donato, P.}
\newblock {\em An Introduction to Homogenization}.
\newblock Oxford University Press, Oxford, 1999.

\bibitem{EEn03}
{\sc E, W., and Engquist, B.}
\newblock The heterogeneous multiscale methods.
\newblock {\em Commun. Math. Sci. 1}, 1 (2003), 87--132.

\bibitem{EEn05}
{\sc E, W., and Engquist, B.}
\newblock The heterogeneous multi-scale method for homogenization problems.
\newblock In {\em Multiscale Methods in Science and Engineering}, vol.~44 of
  {\em Lecture Notes in Computational Science and Engineering}. Springer Berlin
  Heidelberg, 2005, pp.~89--110.

\bibitem{EEnHu03}
{\sc E, W., Engquist, B., and Huang, Z.}
\newblock Heterogeneous multiscale method: A general methodology for multiscale
  modeling.
\newblock {\em Physical Review B 67}, 9 (March 2003), 1--4.

\bibitem{EMiZh05}
{\sc E, W., Ming, P., and Zhang, P.}
\newblock Analysis of the heterogeneous multiscale method for elliptic
  homogenization problems.
\newblock {\em J. Amer. Math. Soc. 18}, 1 (2005), 121--156.

\bibitem{KoWi14}
{\sc Kohn, R.~V., and Wirth, B.}
\newblock Optimal fine-scale structures in compliance minimization for a
  uniaxial load.
\newblock {\em Proceedings of the Royal Society of London A: Mathematical,
  Physical and Engineering Sciences 470}, 2170 (2014).

\bibitem{KoWi15}
{\sc Kohn, R.~V., and Wirth, B.}
\newblock Optimal fine-scale structures in compliance minimization for a shear
  load.
\newblock {\em Communications in Pure and Applied Mathematics\/} (2015).
\newblock to appear.

\bibitem{Lu16}
{\sc L{\"u}then, N.}
\newblock Numerical shape optimization of branching-periodic elastic
  structures.
\newblock Master thesis, Univerisity of Bonn, 2016.

\bibitem{Mi02}
{\sc Milton, G.~W.}
\newblock {\em The Theory of Composites}.
\newblock Cambridge University Press, 2002.

\bibitem{MoMo77}
{\sc Modica, L., and Mortola, S.}
\newblock Un esempio di {$\Gamma ^{-}$}-convergenza.
\newblock {\em Boll. Un. Mat. Ital. B (5) 14}, 1 (1977), 285--299.

\bibitem{Mu09c}
{\sc M{\"u}ller, R.}
\newblock Hierarchical microimaging of bone structure and function.
\newblock {\em Nat. Rev. Rheumatol. 5}, 7 (2009), 373--381.

\bibitem{MuBrDi15}
{\sc M{\"u}ller, V., Brylka, B., Dillenberger, F., Gl{\"o}ckner, R., and
  B{\"o}hlke, T.}
\newblock Homogenization of elastic properties of short-fiber reinforced
  composites based on microstructure data.
\newblock {\em J. Compos. Mater. 50}, 3 (mar 2015), 297--312.

\bibitem{NaWiAm11}
{\sc Nemat-Nasser, S., Wills, J., Srivastava, A., and Amirkhizi, A.}
\newblock Homogenization of periodic elastic composites and locally resonant
  materials.
\newblock {\em Phys. Rev. B 83\/} (Mar 2011), 104103.

\bibitem{PeRuWi12}
{\sc Penzler, P., Rumpf, M., and Wirth, B.}
\newblock A phase-field model for compliance shape optimization in nonlinear
  elasticity.
\newblock {\em ESAIM: Control, Optimisation and Calculus of Variations 18}, 1
  (2012), 229--258.

\end{thebibliography}

\end{document}